\documentclass[11pt,a4paper]{amsart}
\usepackage[latin1]{inputenc}
\usepackage[TS1,T1]{fontenc}
\usepackage[latin1]{inputenc}
\usepackage{amsfonts,amssymb,latexsym}
\usepackage{amsmath,enumerate}
\usepackage{amsthm}
\usepackage{amsfonts,latexsym}
\usepackage[latin1]{inputenc}

\newtheorem{thm}{Theorem}[section]
\newtheorem{coro}{Corollary}[section]

\newtheorem{Step}{Step}
\newtheorem{prop}{Proposition}[section]

\newcommand{\R}{\mathbb{R}}
\newcommand{\N}{\mathbb{N}}
\newcommand{\p}{\partial}

\begin{document}

\title[Sobolev trace constant]{Estimates for the Sobolev trace constant
with critical exponent and applications}

\author[J. Fern\'andez Bonder and N. Saintier]{Juli\'an Fern\'andez Bonder and
Nicolas Saintier}

\address{Departamento de Matem\'atica, FCEyN UBA (1428)
\hfill\break\indent Buenos Aires, Argentina. }

\email{JFB: {\tt jfbonder@dm.uba.ar}, NS: {\tt saintier@math.jussieu.fr} }

\keywords{Sobolev trace embedding, Optimal design problems, Critical exponents.
\\
\indent 2000 {\it Mathematics Subject Classification.} 35J20, 35P30, 49R50}


\begin{abstract}
In this paper we find estimates for the optimal constant in the critical
Sobolev trace inequality $S\|u\|_{L^{p_*}(\partial\Omega)}^p \le
\|u\|_{W^{1,p}(\Omega)}^p$ that are independent of $\Omega$. This estimates
generalized those of \cite{AY} for general $p$. Here $p_* := p(N-1)/(N-p)$ is
the critical exponent for the immersion and $N$ is the space dimension.

Then we apply our results first to prove existence of positive solutions to a
nonlinear elliptic problem with a nonlinear boundary condition with critical
growth on the boundary, generalizing the results of \cite{BR}. Finally, we
study an optimal design problem with critical exponent.
\end{abstract}

\maketitle

\section{Introduction}

Sobolev inequalities are relevant for the study of boundary value problems for
differential operators. They have been studied by many authors and it is by now
a classical subject. It at least goes back to \cite{Aubin1}, for more
references see \cite{DH}. In particular, the Sobolev trace inequality has been
intensively studied in \cite{Biezuner, Escobar, FBLDR, BR, LZ}, etc.

Let $\Omega$ be a bounded smooth domain of $\R^N$. For any $1<p<N$, the Sobolev
trace immersion says that there exists a constant $S>0$ such that
$$
S\Big(\int_{\partial\Omega} |u|^{p_*}\, dS\Big)^{p/p_*} \le \int_{\Omega}
|\nabla u|^p + |u|^p\, dx
$$
for any $u\in W^{1,p}(\Omega)$, where $W^{1,p}(\Omega)$ is the usual Sobolev spaces of the functions $u\in L^p(\Omega$ such that
$\nabla u\in L^p(\Omega$. Here $p_* := p(N-1)/(N-p)$ is the critical
exponent for this inequality.

The optimal constant in the above inequality is the largest possible $S$, that
is
$$
S = S_p(\Omega) := \inf \frac{\displaystyle \int_{\Omega} |\nabla u|^p +
|u|^p\, dx}{\displaystyle \Big(\int_{\partial\Omega} |u|^{p_*}\,
dS\Big)^{p/p_*}},
$$
where the infimum is taken over the set $X := W^{1,p}(\Omega)\setminus
W_0^{1,p}(\Omega)$, $W_0^{1,p}(\Omega)$ being the closure for the $W^{1,p}$-norm of the space of smooth functions with compact support in $\Omega$.

The dependance of $S$ with respect to $p$ and $\Omega$ has been studied by many
authors, specially in the {\em subcritical case}, i.e. where $p_*$ is replaced
by any exponent $q$ such that $1<q<p_*$. See, for instance \cite{FdP, FBFR} and
references therein.

The analysis for the critical case is more involved because the immersion
$W^{1,p}(\Omega)\hookrightarrow L^{p_*}(\partial\Omega)$ is no longer compact
and so the existence of minimizers for $S$ does not follows by standard
methods.

To overcome this problem, in \cite{BR}, the authors use an old idea from T.
Aubin \cite{Aubin1}. In fact, let $K_p^{-1}$ be the best trace constant for the
embedding $W^{1,p}(\R^n_+)\hookrightarrow L^{p_*}(\partial\R^n_+)$, namely
\begin{equation}\label{BestConstant}
K_p^{-1}=\inf_{u\in W^{1,p}(\R^n_+)\setminus W^{1,p}_0(\R^n_+} \frac{ \int_{\R^n_+} |\nabla u|^p dx }{
\left(\int_{\p\R^n_+} |u|^{p_*} dS\right)^{p/p_*} }.
\end{equation}

In \cite{BR} it is shown, following ideas from \cite{Aubin1}, that if
\begin{equation}\label{Cond}
S_p(\Omega)  < K_p^{-1},
\end{equation}
then there exists an extremal for $S_p(\Omega)$. Taking the function $u\equiv
1$ in the definition of $S_p(\Omega)$ one obtain that if
$$
\frac{|\Omega|}{|\p\Omega|^\frac{p}{p_*}}< K_p^{-1},
$$
then \eqref{Cond} is satisfied. Observe that this is a global condition on
$\Omega$.

It follows from Lions \cite{Lions} that the infimum \eqref{BestConstant} is
achieved. The value of $K_p$ is explicitely known when $p=2$ (see Escobar
\cite{Escobar}).

Recently, Biezuner \cite{Biezuner} proved that $K_p$ is also the best first
constant in the inequality,
$$
\left( \int_{\p\Omega} |u|^{p_*} dS \right)^\frac{p}{p_*} \le A\int_\Omega
|\nabla u|^p dx + B\int_\Omega |u|^p dx,
$$
in the sense that, for any $\epsilon>0$, there exists a constant $C_\epsilon$
such that
\begin{equation}\label{OptimalInequ}
\left( \int_{\p\Omega} |u|^{p_*} dS \right)^\frac{p}{p_*} \le
(K_p+\epsilon)\int_\Omega |\nabla u|^p dx + C_\epsilon\int_\Omega |u|^p dx,
\end{equation}
for every  $u\in W^{1,p}(\Omega)$, and $K_p$ is the lowest possible constant.
This fact will be used in a crucial way in the course of the paper.

On the other hand a local condition ensuring \eqref{Cond}, depending only on
local geometric properties of $\Omega$, is known to hold in the case $p=2$.
Indeed Adimurthi-Yadava \cite{AY} obtained \eqref{Cond} assuming the existence
of a ``good point'' $x\in\p\Omega$, i.e. a point $x$ at which the mean
curvature of $\p\Omega$ is positive and such that, in a neighborhood of $x$,
$\Omega$ lies on one side of the tangent plane at $x$. The method in their
proof is the use as test-functions of a suitable rescaling of the extremals of
\eqref{BestConstant}.

These extremals are explicitly known for $p=2$ since Escobar's work
\cite{Escobar} who conjectured the result for any $p\in (1,N)$. This conjecture
has recently been proved by Nazaret \cite{Nazaret} using a mass-transportation
method. It turns out that all the extremals of \eqref{BestConstant} are of the
form
\begin{equation}\label{Extremals}
\begin{aligned}
U_{\epsilon,y_0}(y,t)
 = & \frac{\epsilon^\frac{N-p}{p(p-1)}}{[(t+\epsilon)^2+|y-y_0|^2]^\frac{N-p}{2(p-1)}} \\
 = & \epsilon^{-\frac{N-p}{p}}U\left(\frac{y-y_0}{\epsilon},\frac{t}{\epsilon}\right)
\end{aligned}
\end{equation}
where $\epsilon>0$ and $y, y_0\in\R^{N-1} = \partial\R^N_+$, $t>0$, with
\begin{equation}\label{Def_U}
 U(y,t)= \frac{1}{[(t+1)^2+|y|^2]^\frac{N-p}{2(p-1)}}.
\end{equation}

The knowledge of this extremals allows us first to compute the explicit value
of $K_p$:
\begin{prop}\label{Prop_ValueBestConstant}
The value of $K_p$ is
$$
K_p^{-1} = \left(\frac{N-p}{p-1}\right)^{p-1} \pi^\frac{p-1}{2} \left(
\frac{\Gamma\left(\frac{N-1}{2(p-1)}\right)}
{\Gamma\left(\frac{p(N-1)}{2(p-1)}\right)} \right)^\frac{p-1}{N-1}.
$$
\end{prop}

Applying a similar technique as in \cite{AY}, we can use the rescaled extremals
for $K_p$ and obtain a local (geometrical) condition on $\Omega$ such that
\eqref{Cond} is satisfied.

In fact, we can deal with a slightly more general problem. Namely
\begin{equation}\label{MainPb1}
\lambda = \lambda(p,\Omega) := \inf\frac{\displaystyle \int_\Omega |\nabla u|^p
+ h(x) |u|^p\, dx}{\displaystyle \Big(\int_{\p\Omega} |u|^{p_*}\,
dS\Big)^{p/p_*}}
\end{equation}
where the infimum is taken over $X$ and the function $h\in
C^1(\overline{\Omega})$ is such that there exists $c>0$ satisfying
\begin{equation}\label{Coercivity}
\int_\Omega |\nabla u|^p + h(x) |u|^p\, dx \ge c \|u\|_{W^{1,p}(\Omega)}^p
\end{equation}
for any $u\in X$.

We are lead to the following generalization of the notion of ``good point'' to
our case: we say that a point $x\in\partial\Omega$ is a ``good point'' if there
exists $r>0$ such that $\Omega\cap B_r(x)$ lies on one side of the tangent
plane at $x$ and either $H(x)>0$ or, if $H(x)=0$, either
$$
h(x)<0 \text{ if } N=2,3,4 \text{ and } p<\sqrt{N}
$$
or, if $N\ge 5$,
\begin{equation*}
\begin{aligned}
& h(x)<0 \text{ if } p<2, \\
& \frac{N}{N-1}\sum \lambda_i^2 - 2\sum_{i<j}\lambda_i\lambda_j <
\frac{-8(N-1)h(x)}{(N-2)(N-4)} \text{ if } p=2, \\
& \frac{p+N-2}{N-1}\sum \lambda_i^2 - 2\sum_{i<j}\lambda_i\lambda_j<0 \text{ if
} 2<p<(N+2)/3.
\end{aligned}
\end{equation*}
where the $\lambda_i$'s are the principal curvatures at $x$ and $H(x)$ is the
mean curvature at $x$.

Remark that our method gives the restriction $1<p<(N+1)/2$ and also that a
``good point'' in the sense of Adimurthi-Yadava is also a ``good point'' in our
sense.

We get the following theorem:

\begin{thm}\label{thm1}
Let $1<p<(N+1)/2$. If there exist a ``good point'' $x\in\partial\Omega$,
then
\begin{equation}\label{eq-thm1}
\lambda < K_p^{-1}.
\end{equation}
\end{thm}

As a consequence of Theorem \ref{thm1} we have
\begin{coro}\label{coro-thm1}
Under the hypotheses of Theorem \ref{thm1}, the infimum \eqref{MainPb1} is
achieved.
\end{coro}

Observe that any extremal $u$ can be taken to be nonnegative (just replace $u$
by $|u|$), and if we take it {\em normalized} as
$\|u\|_{L^{p_*}(\partial\Omega)}=1$, it is an eigenfunction associated to the
eigenvalue $\lambda$ in the sense that it is a weak solution of the following
Steklov-like eigenvalue problem
\begin{equation}\label{MainEqu}
\begin{cases}
-\Delta_pu + h(x) u^{p-1} = 0 \quad &\text{in } \Omega \\
|\nabla u|^{p-2} \frac{\p u}{\p\nu} = \lambda u^{p_*-1} &\text{on } \p\Omega
\end{cases}
\end{equation}
where $\Delta_pu = \text{div}(|\nabla u|^{p-2}\nabla u)$ is the $p-$Laplacian
and $\nu$ is the unit outward normal of $\Omega$.

Then it follows by the results of Cherrier \cite{Cherrier} that $u$ is smooth
on $\Omega$ and continuous up to the boundary. Moreover, it is strictly
positive in $\overline{\Omega}$ (see, for instance, \cite{FBR-JMAA}) so any
extremal has constant sign.

As an application of Theorem \ref{thm1}, we study a shape optimization problem
related to $\lambda$. Given $\alpha\in (0,|\Omega|)$, where $|\Omega|$ denotes
the volume of $\Omega$, and a measurable subset $A\subset\Omega$ of volume
$\alpha$, we first consider the minimization problem
\begin{equation}\label{MainPbShape1}
 \lambda_A = \inf\frac{\displaystyle \int_\Omega |\nabla u|^p + h(x) |u|^p\,
dx}{\displaystyle \Big(\int_{\p\Omega} |u|^{p_*}\, dS\Big)^{p/p_*}}
\end{equation}
where the infimum is taken over $X_A := \{u\in X\ |\ u|_A = 0 \text{ a.e.}\}$
and the function $h\in C^1(\overline{\Omega})$ is such that the coercivity
assumption \eqref{Coercivity} holds

As a consequence of Theorem \ref{thm1}, we have

\begin{thm}\label{thmShape1}
Let $1<p<(N+1)/2$ and let $A\subset \Omega$ be such that $|A|=\alpha$. Assume
that there exists a ``good point'' $x\in\p\Omega$ such that $B_r(x)\cap
A=\emptyset$ for some $r>0$. Then $\lambda_A$ is attained by some nonnegative
nontrivial $u_A$.
\end{thm}

These extremals $u_A$ are eigenfunctions associated to the eigenvalue
$\lambda_A$ in the sense that, if $A$ is closed, they are weak solutions of the
following Steklov-like eigenvalue problem
\begin{equation}\label{MainEquShape}
 \begin{cases}
  -\Delta_p u + h(x) u^{p-1} = 0 & \text{in } \Omega\setminus A \\
  |\nabla u|^{p-2} \frac{\p u}{\p\nu} = \lambda_A u^{p_*-1} & \text{on } \p\Omega\setminus A \\
  u = 0 & \text{in } A
 \end{cases}
\end{equation}

We consider the following shape optimization problem:

\begin{itemize}
\item[]{\em For a fixed $0<\alpha<|\Omega|$, find a set $A_*$ of measure
$\alpha$ that minimizes $\lambda_A$ among all measurable subsets $A\subset
\Omega$ of measure $\alpha$. That is, }
$$
\lambda(\alpha) := \inf_{A\subset\Omega, |A|=\alpha}{\lambda_A} =
\lambda_{A_*}.
$$
\end{itemize}

In this paper we prove that there exist an optimal set $A_*$ (with their
corresponding extremals $u_*$) for this optimization problem.

This optimization problem in the subcritical case (that is, when $p_*$ is
replaced by an exponent $q$ with $1<q<p_*$) has been considered recently. In
fact, in \cite{BRW1} the existence of an optimal set has been established, see
also \cite{BGR} for numerical computations. Then, in \cite{BRW2}, the interior
regularity of optimal sets was analyzed in the case $p=2$. We remark that in
the result of \cite{BRW2} the subcriticality plays no role, so this local
regularity result holds true also for this critical case.

We prove,
\begin{thm}\label{thmShape2}
Let $1<p<(N+1)/2$. If there exists a ``good point'' $x\in\partial\Omega$, then
$\lambda(\alpha)$ is achieved.
\end{thm}

Problems of optimal design related to eigenvalue problems like
\eqref{MainEquShape} appear in several branches of applied mathematics,
specially in the case $p=2$. For example in problems of minimization of the
energy stored in the design under a prescribed loading. We refer to \cite{CC}
for more details.

We want to stress that Theorem \ref{thmShape2} is new, even in the case $p=2$.

\subsection*{Organization of the paper}
In the next section we deal with the proof of the applications of the estimate
$\lambda < K_p^{-1}$, that is, we deal with the proof of Corollary
\ref{coro-thm1} and Theorems \ref{thmShape1} and \ref{thmShape2}. We leave for
the final section the computation of $K_p$ and the proof of Theorem \ref{thm1}.

\section{Applications of Theorem \ref{thm1}}
\setcounter{equation}{0}

In this section we use Theorem \ref{thm1}, that is proved in the Section 3, and
prove Corollary \ref{coro-thm1}, Theorem \ref{thmShape1} and Theorem
\ref{thmShape2}.

\subsection{Proof of Corollary \ref{coro-thm1}}

We first prove that $\lambda$ is attained as soon as \eqref{eq-thm1} is
satisfied. Since this kind of criterion is classical (see e.g. \cite{DN} or
\cite{BR}), we only sketch the proof for the reader's convenience.

Let $\{u_n\}_{n\in\N}\subset X$ be a minimizing sequence for (\ref{MainPb1})
normalized such that $\|u_n\|_{L^{p_*}(\partial\Omega)}=1$. According to
(\ref{Coercivity}), this sequence is bounded in $X$ and thus it converges up to
a subsequence to some $u\in X$ weakly in $X$, strongly in $L^p(\Omega)$ and
a.e.

Using Ekeland's variational principle (see \cite{Willem} Theorems 8.5 and
8.14), we can assume that $\{u_n\}_{n\in\N}$ is a Palais-Smale sequence for the
functional $J:W^{1,p}(\Omega)\to \R$ defined by
$$
J(u)= \frac{1}{p} \int_\Omega |\nabla u|^p + h(x)|u|^p\, dx -
\frac{\lambda}{p_*} \int_{\p\Omega} |u|^{p_*}\, dS,
$$
in the sense that the sequence $\{J(u_n)\}_{n\in\N}$ is bounded and $DJ(u_n)\to
0$ strongly in $(W^{1,p}(\Omega))^*$. Letting $v_n:=u_n-u$, we can also assume
that, up to a subsequence,
$$
|v_n|^{p_*}\, dS \rightharpoonup d\nu,\qquad |\nabla v_n|^p\, dx
\rightharpoonup d\mu,
$$
weakly in the sense of measures, where $\mu$ and $\nu$ are nonnegative measures
such that $\text{supp}(\nu)\subset\p\Omega$.

According to \eqref{OptimalInequ}, we have for any $\phi\in
C^1(\overline{\Omega})$ that
$$
\left( \int_{\p\Omega} |\phi v_n|^{p_*}\, dS \right)^{p/p_*} \le
(K_p+\epsilon)\int_\Omega |\nabla (\phi v_n)|^p\, dx + C_\epsilon \int_\Omega
|\phi v_n|^p\, dx.
$$
Passing to the limit in this expression, first in $n\to\infty$ and then in
$\epsilon\to 0$, we get that
$$
\left( \int_{\p\Omega} |\phi|^{p_*}\, d\nu \right)^{p/p_*} \le K_p \int_\Omega
|\phi|^p\, d\mu
$$
for any $\phi\in C^1(\overline{\Omega})$. From this inequality, we can deduce
as in \cite{Lions} Lemma 2.3, the existence of a sequence of points
$\{x_i\}_{i\in I}\subset\p\Omega$, $I\subset\N$, and two sequences of positive
real numbers $\{\nu_i\}_{i\in I}$,  $\{\mu_i\}_{i\in I}$ such that
$$
\nu=\sum_{i\in I} \nu_i\delta_{x_i},\quad \mu\ge\sum_{i\in I} \mu_i\delta_{x_i}
\quad \text{and}\quad \mu_i\ge K_p^{-1}\nu_i^{p/p_*} \quad \forall\ i\in I.
$$
Therefore,
\begin{equation}\label{CCP}
\begin{cases}
|u_n|^{p_*}dS & \rightharpoonup  |u|^{p_*}dS + \sum_{i\in I} \nu_i\delta_{x_i} \\
|\nabla u_n|^pdx & \rightharpoonup  |\nabla u_n|^pdx + \mu \ge |\nabla u_n|^pdx + \sum_{i\in I} \mu_i\delta_{x_i} \\
\mu_i & \ge K_p^{-1}\nu_i^{p/p_*}~~\forall~i\in I.
\end{cases}
\end{equation}
It can also be shown that $\{v_n\}_{n\in\N}$ is a Palais-Smale sequence for the
functional $I:W^{1,p}(\Omega)\to \R$ defined by
$$
I(u):=J(u)-\int_\Omega h(x)|u|^p\, dx
$$
(see e.g. \cite{Saintier}). In particular, for any $\phi\in
C^1(\overline{\Omega})$,
\begin{align*}
 o(1) & = DI(v_n)(v_n\phi) \\
      & = \int_\Omega |\nabla v_n|^{p-2} \nabla v_n\nabla (v_n\phi)\, dx -
      \lambda \int_{\p\Omega} |v_n|^{p_*}\phi\, dS.
\end{align*}
Passing to the limit, we get that $\int_\Omega \phi\, d\mu =
\lambda\int_{\p\Omega} \phi\, d\nu$ for any $\phi\in C^1(\overline{\Omega})$.
Hence $\mu = \lambda\nu$. Using \eqref{CCP}, we then obtain the estimates
\begin{equation}\label{Estimates_mu}
\nu_i\ge (\lambda K_p)^{-\frac{n-1}{p-1}}, \quad \mu_i\ge K_p^{-1}(\lambda
K_p)^{-\frac{n-1}{p-1}}\quad \forall\ i\in I.
\end{equation}
Now, by \eqref{CCP}, \eqref{Coercivity} and \eqref{Estimates_mu}, we arrive at
\begin{align*}
\lambda & = \int_\Omega |\nabla u_n|^p\, dx + \int_\Omega h(x)|u_n|^p\, dx +
o(1) \ge \sum_{i\in I} \mu_i \\
& \ge card(I) K_p^{-1}(\lambda K_p)^{-\frac{n-1}{p-1}}.
\end{align*}
We deduce that if \eqref{eq-thm1} holds, then $I$ is empty. In that case,
$u_n\to u$ strongly in $W^{1,p}(\Omega)$ and in $L^{p_*}(\p\Omega)$. In
particular $u$ is a minimizer for $\lambda$.

This completes the proof \qed

\subsection{Proof of Theorem \ref{thmShape1}}

Arguing exactly as in the proof of Theorem \ref{thm1} we obtain that a
normalized minimizing sequence $\{u_n\}_{n\in\N}\subset X_A$ for $\lambda_A$
converges, up to a subsequence, strongly in $W^{1,p}(\Omega)$ to some $u_A$ as
soon as
\begin{equation}\label{Cond2}
\inf_{u\in X_A} \frac{\displaystyle \int_{\Omega} |\nabla u|^p + |u|^p\,
dx}{\displaystyle \Big(\int_{\partial\Omega} |u|^{p_*}\, dS\Big)^{p/p_*}} <
K_p^{-1}.
\end{equation}
Since there exists a ``good point'' $x\in\p\Omega$ such that $B_r(x)\cap
A=\emptyset$, we deduce from the computations in the next section, by choosing
a cut-off function $\phi$ with support in $B_{r/2}(x)$ in the definition of the
test function $u_\epsilon$ \eqref{test.function}, that this strict inequality
\eqref{Cond2} holds. Hence $u_n\to u$ strongly in $W^{1,p}(\Omega)$ and
$L^{p_*}(\p\Omega)$ and also a.e.. In particular $u$ is a minimizer for $\lambda_A$. \qed

\subsection{Proof of Theorem \ref{thmShape2}}

We begin by noticing that
$$
\lambda(\alpha) = \inf \{\lambda_A,~A\subset\Omega~\text{measurable},~|A|\ge\alpha\}.
$$
Hence
$$
\lambda(\alpha) = \inf_{u\in X,~|\{u=0\}|\ge \alpha} \frac{\displaystyle
\int_{\Omega} |\nabla u|^p + |u|^p\, dx}{\displaystyle
\Big(\int_{\partial\Omega} |u|^{p_*}\, dS\Big)^{p/p_*}}.
$$
Since $\alpha<|\Omega|$ and there exists a ``good point'', it follows from the
test functions computations of the next section, by choosing a function $\phi$
with support in a ball of radius small enough in the definition of $u_\epsilon$
\eqref{test.function}, that $\lambda(\alpha)<K_p^{-1}$.

By the same argument as before, this implies the existence of a nonnegative
$u_*\in X$, $|\{u_*=0\}|\ge \alpha$, such that
$$
\frac{\displaystyle \int_{\Omega} |\nabla u_*|^p + |u_*|^p\, dx}{\displaystyle
\Big(\int_{\partial\Omega} |u_*|^{p_*}\, dS\Big)^{p/p_*}} = \lambda(\alpha).
$$
We now conclude as in \cite{BRW1}, Theorem 1.2, that in fact
$|\{u_*=0\}|=\alpha$ and so $A_* = \{u_*=0\}$ is an optimal set for
$\lambda(\alpha)$.\qed

\section{Proof of Theorem \ref{thm1}}
\setcounter{equation}{0}

In this section we prove our main result. First we recall some very well known
formulae and prove Proposition \ref{Prop_ValueBestConstant}. Finally we prove
Theorem \ref{thm1}.

In all the subsequent computations, the following well known formulae will be
used frequently:
\begin{align*}
& \omega_{N-1} = \text{volume of the standard unit sphere $S^{N-1}$ of $\R^N$}
= \frac{2\pi^\frac{N}{2}}{\Gamma\left(\frac{N}{2}\right)}, \\
\\
& \int_0^{+\infty} \frac{r^\alpha}{(1+r^2)^\beta} dr =
\frac{\Gamma\left(\frac{\alpha+1}{2}\right)\Gamma\left(\frac{2\beta - \alpha-1}{2}\right)}
{2\Gamma(\beta)} \quad \text{ for } 2\beta-\alpha>1,\\
\\
& \Gamma(z)\Gamma(z+\frac{1}{2}) = 2^{1-2z}\sqrt{\pi}\Gamma(2z) \quad \text{
for } Re(z)>0.
\end{align*}

We first compute the value of $K_p$:

\begin{proof}[{\bf Proof of Propostion \ref{Prop_ValueBestConstant}}]
Let $U$ be the function defined by \eqref{Def_U}. We first compute the
$L^{p_*}$-norm of $U$ restricted to $\R^{N-1}\times\{0\} = \p\R^N_+$.
\begin{align*}
\int_{\R^{N-1}} |U(y,0)|^{p_*}\, dy & =  \int_{\R^{N-1}}
\frac{dy}{(1+|y|^2)^{p(N-1)/2(p-1)}}\\
& = \omega_{N-2} \int_0^{\infty} \frac{r^{N-2}\, dr}{(1+r^2)^{p(N-1)/2(p-1)}} \\
 & = \pi^{(N-1)/2} \frac{\Gamma\left(\frac{N-1}{2(p-1)}\right)}{\Gamma\left(\frac{p(N-1)}{2(p-1)}\right)}
\end{align*}
We now compute the $L^p$-norm of the gradient of $U$. First
$$
\nabla U(y,t)= -\frac{N-p}{p-1}
\frac{(y,t+1)}{[(1+t)^2+|y|^2]^{\frac{N-p}{2(p-1)}+1}}.
$$
Using the change of variable $y=(1+t)z$ and passing to polar coordinates, we
can then write
\begin{align*}
\int_{\R^N_+} |\nabla U(y,t)|^p\, dydt
 & = \left(\frac{N-p}{p-1}\right)^p  \int_{\R^N_+} \frac{dydt}{ [(1+t)^2+|y|^2]^{\frac{p(N-1)}{2(p-1)}} } \\
 & = \left(\frac{N-p}{p-1}\right)^p \int_0^{+\infty} \frac{dt}{(1+t)^\frac{N-1}{p-1}} \omega_{N-2}
      \int_0^{+\infty} \frac{r^{N-2}\, dr}{(1+r^2)^\frac{p(N-1)}{2(p-1)}} \\
 & = \left(\frac{N-p}{p-1}\right)^{p-1} \pi^\frac{N-1}{2}
     \frac{\Gamma\left(\frac{N-1}{2(p-1)}\right)}{\Gamma\left(\frac{p(N-1)}{2(p-1)}\right)}.
\end{align*}
Hence
$$
K_p^{-1} = \frac{\displaystyle \int_{\R^N_+} |\nabla U(y,t)|\, dydt}
{\displaystyle \Big(\int_{\R^{N-1}} |U(y,0)|^{p_*} dy \Big)^\frac{p}{p_*}} =
\left(\frac{N-p}{p-1}\right)^{p-1} \pi^\frac{p-1}{2}
\left(\frac{\Gamma\left(\frac{N-1}{2(p-1)}\right)}
{\Gamma\left(\frac{p(N-1)}{2(p-1)}\right)} \right)^\frac{p-1}{N-1}
$$
and the proof is complete
\end{proof}

We now turn our attention to the proof of Theorem \ref{thm1}. Let
$x_0\in\p\Omega$ be a ``good point''. By taking an appropriate chart, we can
assume that $x_0=0$ and that there exist $r>0$ and
$\lambda_1,\dots,\lambda_{N-1}\in\R$ such that
\begin{align*}
B_r\cap\Omega = & \{(y,t)\in B_r,\ t>\rho(y)\} \\
B_r\cap\p\Omega =&  \{(y,t)\in B_r,\ t=\rho(y)\}
\end{align*}
where $y=(y_1,\dots,y_{N-1})\in\R^{N-1}$, $B_r$ is the Euclidean ball centered
at the origin and of radius $r$, and
$$
\rho(y)=\frac{1}{2}\sum_{i=1}^{N-1} \lambda_i y_i^2 + \sum_{i,j,k} c_{ijk}y_i
y_j y_k + O(|y|^4).
$$
Since $x_0=0$ is a ``good point'', we have $\rho\ge 0$. Moreover, the
$\lambda_i$'s are the principal curvatures at $0$ and thus
$$
H(0)=\frac{1}{N-1}\sum_{i=1}^{N-1} \lambda_i.
$$
Let $\phi$ be a smooth radial function with compact support in $B_{r/2}$ be
such that $\phi\equiv 1$ in $B_{r/4}$. We consider the test functions
\begin{equation}\label{test.function}
u_\epsilon(y,t)=\frac{\phi(y,t)}{[(t+\epsilon)^2+|y|^2]^\frac{N-p}{2(p-1)}},~~~\epsilon>0.
\end{equation}
In order to give the asymptotic development of the Rayleigh quotient for
$u_\epsilon$, we first compute the different terms involved:

\begin{Step}\label{step_calcul1}
We have the following estimates:
\begin{equation}\label{NormeGradient}
\int_\Omega |\nabla u_\epsilon|^p\, dx =  A_1 \epsilon^{-\frac{N-p}{p-1}}
        +  \begin{cases}
           A_2 \epsilon^{1-\frac{N-p}{p-1}} + A_3 \epsilon^{2-\frac{N-p}{p-1}} \\
           \hspace{1cm} +\begin{cases} O(\epsilon^{3-\frac{N-p}{p-1}}) \text{ if } p<\frac{N+3}{4} \\
                                        O(\ln(1/\epsilon)) \text{ if }p=\frac{N+3}{4} \\
                                        O(1) \text{ if } \frac{N+3}{4}<p<\frac{N+1}{2}
                          \end{cases} \\
          A_2'\ln(1/\epsilon) \text{ if } p=\frac{N+1}{2} \\
          O(1) \text{ if } p>\frac{N+1}{2}
         \end{cases}
\end{equation}

\begin{equation}\label{NormeLp}
\int_\Omega h(x) |u_\epsilon|^p\, dx = \begin{cases}
       D\epsilon^{-\frac{N-p^2}{p-1}} + \begin{cases}
                                          O(\epsilon^{1-\frac{N-p^2}{p-1}}) \text{ if }p<\frac{-1+\sqrt{4N+5}}{2} \\
                                          O(\ln(1/\epsilon)) \text{ if } p=\frac{-1+\sqrt{4N+5}}{2} \\
                                          O(1) \text{ if } \sqrt{N}>p>\frac{-1+\sqrt{4N+5}}{2}
                                         \end{cases} \\
       O(\ln(1/\epsilon)) \text{ if } p=\sqrt{N} \\
       O(1) \text{ if } p>\sqrt{N}
      \end{cases}
\end{equation}

\begin{equation}\label{NormeCritique}
\begin{aligned}
\int_{\p\Omega} |u_\epsilon|^{p_*}\, dS =&
B_1 \epsilon^{-1-\frac{N-p}{p-1}} + B_2 \epsilon^{-\frac{N-p}{p-1}}   \\
& + \begin{cases}
           B_3 \epsilon^{1-\frac{N-p}{p-1}} +
           \begin{cases}
            O(\epsilon^{2-\frac{N-p}{p-1}}) \text{ if } p<\frac{N+2}{3} \\
            O(\ln(1/\epsilon)) \text{ if } p=\frac{N+2}{3} \\
            O(1) \text{ if } \frac{N+2}{3}<p<\frac{N+1}{2}
           \end{cases} \\
          B_4 \ln(1/\epsilon) \text{ if } p=\frac{N+1}{2} \\
          O(1) \text{ if } p>\frac{N+1}{2}
         \end{cases}
\end{aligned}
\end{equation}
where
$$
A_1 = \frac{1}{2}\left(\frac{N-p}{p-1}\right)^{p-1}\omega_{N-2}
\frac{\Gamma\left(\frac{N-1}{2}\right)\Gamma\left(\frac{N-1}{2(p-1)}\right)}
{\Gamma\left(\frac{p(N-1)}{2(p-1)}\right)}
$$
$$
A_2 = -\frac{H(0)\omega_{N-2}}{4} \left(\frac{N-p}{p-1}\right)^p
\frac{\Gamma\left(\frac{N+1}{2}\right)\Gamma\left(\frac{N-2p+1}{2(p-1)}\right)}
{\Gamma\left(\frac{p(N-1)}{2(p-1)}\right)}
$$
$$
A_2' = -\frac{H(0)\omega_{N-2}}{2} \left(\frac{N-p}{p-1}\right)^p
$$
$$
A_3 = \frac{\omega_{N-2}}{16} \left(\frac{N-p}{p-1}\right)^p
\frac{\Gamma\left(\frac{N-1}{2}\right)\Gamma\left(\frac{N-2p+1}{2(p-1)}\right)}
{\Gamma\left(\frac{p(N-1)}{2(p-1)}\right)}
\left(\frac{3}{2}\sum\lambda_i^2+\sum_{i<j}\lambda_i\lambda_j \right)
$$
$$
B_1 = \omega_{N-2} \frac{\Gamma\left(\frac{N-1}{2}\right)
\Gamma\left(\frac{N-1}{2(p-1)}\right)}{2\Gamma\left(\frac{p(N-1)}{2(p-1)}\right)}
$$
$$
B_2 = -\frac{\omega_{N-2} \sum \lambda_i}{8} \frac{p(N-1)}{p-1}
\frac{\Gamma\left(\frac{N-1}{2}\right)\Gamma\left(\frac{N-1}{2(p-1)}\right)}
{\Gamma\left(1+\frac{p(N-1)}{2(p-1)}\right)}
$$
\begin{align*}
B_3 = & \frac{\omega_{N-2}}{32}
\frac{\Gamma\left(\frac{N-1}{2}\right)\Gamma\left(\frac{N-2p+1}{2(p-1)}\right)}
{\Gamma\left(\frac{p(N-1)}{2(p-1)}\right)} \times \\
&\left\{\left(1+\frac{3(N-2p+1)}{p-1}\right) \sum \lambda_i^2 +
\left(-2+\frac{2(N-2p+1)}{p-1}\right) \sum_{i<j} \lambda_i\lambda_j \right\}
\end{align*}
$$
B_4 = \frac{\omega_{N-2}}{2} \left\{\left(\frac{1}{N-1} -
\frac{p(N-1)}{4(p-1)}\right) \sum \lambda_i^2 - \frac{p(N-1)}
{2(p-1)}\sum_{i<j} \lambda_i\lambda_j + o(1) \right\}
$$
$$
D = h(0) \frac{p-1}{N-p^2} \omega_{N-2}
\frac{\Gamma\left(\frac{N-1}{2}\right)\Gamma\left(\frac{N-p^2+p-1}{2(p-1)}\right)}
{2\Gamma\left(\frac{p(N-p)}{2(p-1)}\right)}
$$
\end{Step}

\begin{proof}[Proof of Step \ref{step_calcul1}]
We have
$$
[(t+\epsilon)^2+|y|^2]^\frac{N-1}{p-1} |\nabla u_\epsilon|^2 =
\left(\frac{N-p}{p-1}\right)^2 \phi^2 + |\nabla\phi|^2 -
2\frac{N-p}{p-1}\phi(y\cdot \nabla_y\phi + (t+\epsilon)\p_t\phi)
$$
Hence in $B_{r/4}$,
$$
|\nabla u_\epsilon|^p = \left(\frac{N-p}{p-1}\right)^p
\frac{1}{[(t+\epsilon)^2+|y|^2]^\frac{p(N-1)}{2(p-1)}},
$$
and then
$$
\int_\Omega |\nabla u_\epsilon|^p\, dx =  \left(\frac{N-p}{p-1}\right)^p (I_1 -
I_2) + O(1)
$$
with
$$
I_1 = \int_{Q_a} \frac{1}{[(t+\epsilon)^2+|x|^2]^\frac{p(n-1)}{2(p-1)}}
\quad\text{ and }\quad I_2 = \int_{Q_a\setminus \Omega}
\frac{1}{[(t+\epsilon)^2+|x|^2]^\frac{p(n-1)}{2(p-1)}},
$$
where $Q_a := \{(y,t)\ |\ |y|\le a \text{ and } 0\le t\le a\}$.

Changing variables $y=(1+t)z$ and passing to polar coordinates, we have
\begin{align*}
I_1 & = \int_{Q_a} \frac{1}{[(t+\epsilon)^2+|y|^2]^\frac{p(N-1)}{2(p-1)}}\, dydt  \\
& = \epsilon^{-\frac{N-p}{p-1}} \int_{\R^N_+}  \frac{1}{[(1+t)^2+|y|^2]^\frac{p(N-1)}{2(p-1)}}\, dydt + O(1) \\
& = \epsilon^{-\frac{N-p}{p-1}} \omega_{N-2} \int_0^\infty
\frac{dt}{(1+t)^\frac{N-1}{p-1}} \int_0^\infty \frac{r^{N-2}\, dr}
{(1+r^2)^\frac{p(N-1)}{2(p-1)}} + O(1)
\end{align*}
Hence
\begin{equation}\label{I_1}
I_1 = \epsilon^{-\frac{N-p}{p-1}} \frac{p-1}{N-p} \omega_{N-2}
\frac{\Gamma\left(\frac{N-1}{2}\right) \Gamma\left(\frac{N-1}{2(p-1)}\right)}
{2\Gamma\left(\frac{p(N-1)}{2(p-1)}\right)} + O(1).
\end{equation}

On the other hand, according to Taylor's formula,
\begin{align*}
I_2  = & \int_{|y|\le a} \int_0^{\rho(y)}\frac{1}{[(t+\epsilon)^2+|y|^2]^\frac{p(N-1)}{2(p-1)}}\, dtdy \\
= & \int_{|y|\le a} \frac{\rho(y)\, dy}
{(\epsilon^2+|y|^2)^\frac{p(N-1)}{2(p-1)}} - \frac{p(N-1)}{2(p-1)} \epsilon
\int_{|y|\le a} \frac{\rho(y)^2\, dy}{(\epsilon^2+|y|^2)^{\frac{p(N-1)}{2(p-1)}+1}} \\
 & + O\left(\int_{|y|\le a} \frac{|y|^6\, dy}{(\epsilon^2+|y|^2)^{\frac{p(N-1)}{2(p-1)}+1}} \right) \\
 = &\ I_3 - \frac{p(N-1)}{2(p-1)} \epsilon I_4
       + \begin{cases} O\left( \epsilon^{3-\frac{N-p}{p-1}}\right), \text{ if } p<\frac{N+3}{4}  \\
                       O(\ln(1/\epsilon)), \text{ if } p=\frac{N+3}{4}  \\
                       O(1), \text{ if } p>\frac{N+3}{4}
          \end{cases}
\end{align*}
As the sphere is symmetric, we have
\begin{align*}
I_3 & = \frac{1}{2}H(0) \int_{|y|\le a} \frac{|y|^2\, dy}
{(\epsilon^2+|y|^2)^\frac{p(N-1)}{2(p-1)}} + O\left(\int_{|y|\le a}
\frac{|y|^4\, dy}{(\epsilon^2+|y|^2)^\frac{p(N-1)}{2(p-1)}} \right)
\end{align*}
with
\begin{equation}\label{E1}
\begin{aligned}
\int_{|y|\le a} & \frac{|y|^2\, dy}{(\epsilon^2+|y|^2)^\frac{p(N-1)}{2(p-1)}}
 = \epsilon^{1-\frac{N-p}{p-1}} \omega_{N-2} \int_0^{a/\epsilon} \frac{r^N dr}{(1+r^2)^{\frac{p(N-1)}{2(p-1)}}} \\
 & = \begin{cases}
       \epsilon^{1-\frac{N-p}{p-1}}\omega_{N-2}
       \frac{\Gamma\left(\frac{N+1}{2}\right) \Gamma\left(\frac{N-2p+1}{2(p-1)}\right) }
        { 2\Gamma\left(\frac{p(N-1)}{2(p-1)}\right) } + O(1) \text{ if } p<\frac{N+1}{2} \\
       \approx\omega_{N-2}\ln(1/\epsilon) \text{ if } p<\frac{N+1}{2} \\
       O(1) \text{ if } p>\frac{N+1}{2} \\
     \end{cases}
\end{aligned}
\end{equation}
and
\begin{equation}\label{E2}
\begin{aligned}
\int_{|y|\le a} \frac{|y|^4\, dy}{(\epsilon^2+|y|^2)^\frac{p(N-1)}{2(p-1)}}
& = \epsilon^{3-\frac{N-p}{p-1}} \omega_{N-2} \int_0^{a/\epsilon} \frac{r^{N+2}\, dr}
{(1+r^2)^{\frac{p(N-1)}{2(p-1)}}} \\
& = \begin{cases}
       O(\epsilon^{3-\frac{N-p}{p-1}}) \text{ if } p<\frac{N+3}{4} \\
       O(\ln(1/\epsilon)) \text{ if } p=\frac{N+3}{4} \\
       O(1) \text{ if } p>\frac{N+3}{4}
     \end{cases}
\end{aligned}
\end{equation}
Since $\frac{N+3}{4}<\frac{N+1}{2}$ we get
\begin{equation*}
\begin{aligned}
I_3 & = \begin{cases} \epsilon^{1-\frac{N-p}{p-1}}\omega_{N-2} H(0)
\frac{\Gamma\left(\frac{N+1}{2}\right)
\Gamma\left(\frac{N-2p+1}{2(p-1)}\right)}{4\Gamma\left(\frac{p(N-1)}{2(p-1)}\right)}
+ \begin{cases}
           O(\epsilon^{3-\frac{N-p}{p-1}}) \text{ if } p<\frac{N+3}{4} \\
           O(\ln(1/\epsilon)) \text{ if } p=\frac{N+3}{4} \\
           O(1) \text{ if } \frac{N+3}{4}<p<\frac{N+1}{2}
          \end{cases} \\
     \approx\frac{1}{2}H(0)\omega_{N-2}\ln(1/\epsilon) \text{ if } p=\frac{N+1}{2} \\
     O(1) \text{ if } p>\frac{N+1}{2}
   \end{cases}
\end{aligned}
\end{equation*}
Concerning $I_4$, we have
\begin{align*}
I_4 = & \frac{1}{4}\sum \lambda_i^2  \int_{|y|\le a} \frac{y_i^4\, dy}
{(\epsilon^2+|y|^2)^{\frac{p(N-1)}{2(p-1)}+1}}\\
& + \frac{1}{2}\sum_{i<j} \lambda_i\lambda_j  \int_{|y|\le a} \frac{y_i^2
y_j^2\, dy} {(\epsilon^2+|y|^2)^{\frac{p(N-1)}{2(p-1)}+1}} \\
& \qquad + O\left(\int_{|y|\le a} \frac{|y|^5\, dy}
{(\epsilon^2+|y|^2)^{\frac{p(N-1)}{2(p-1)}+1}} \right).
\end{align*}
First we compute
\begin{align*}
\int_{|y|\le a} \frac{y_i^4\, dy}{(\epsilon^2+|y|^2)^{\frac{p(N-1)}{2(p-1)}+1}}
& = \epsilon^{1-\frac{N-p}{p-1}} \int_{|y|\le a/\epsilon} \frac{y_i^4\, dy}
{(\epsilon^2+|y|^2)^{\frac{p(N-1)}{2(p-1)}+1}} \\
& = \begin{cases}
      O(1) \text{ if } p>\frac{N+1}{2} \\
      \approx \omega_{N-2}\ln(1/\epsilon) \text{ if } p=\frac{N+1}{2}
     \end{cases}
\end{align*}
and if $p<\frac{N+1}{2}$,
\begin{align*}
& \int_{|y|\le a} \frac{y_i^4\, dy}{(\epsilon^2+|y|^2)^{\frac{p(N-1)}{2(p-1)}+1}} \\
& = 2 \epsilon^{1-\frac{N-p}{p-1}} \omega_{N-3} \int_0^\infty \frac{r^{N-3}\,
dr} {(1+r^2)^{\frac{p(N-1)}{2(p-1)}-\frac{3}{2}}} \int_0^\infty \frac{y^4\, dy}
{(1+y^2)^{\frac{p(N-1)}{2(p-1)}+1}} + O(1).
\end{align*}

Hence
\begin{equation}\label{I4_1}
\begin{aligned}
& \int_{|y|\le a} \frac{y_i^4\, dy}{(\epsilon^2+|y|^2)^{\frac{p(N-1)}{2(p-1)}+1}} \\
& = \begin{cases}
      \epsilon^{1-\frac{N-p}{p-1}} \frac{\omega_{N-3}}{2}
      \frac{\Gamma\left(\frac{N-2}{2}\right)\Gamma\left(\frac{N-2p+1}{2(p-1)}\right)\Gamma\left(\frac{5}{2}\right)}
           {\Gamma\left(\frac{p(N-1)}{2(p-1)}+1\right)} + O(1) \text{ if } p<\frac{N+1}{2} \\
       \approx \omega_{N-2}\ln(1/\epsilon) \text{ if } p=\frac{N+1}{2} \\
      O(1) \text{ if } p>\frac{N+1}{2}
     \end{cases}
\end{aligned}
\end{equation}
In the same way
$$
\int_{|y|\le a} \frac{y_i^2y_j^2\, dy}
{(\epsilon^2+|y|^2)^{\frac{p(N-1)}{2(p-1)}+1}} = \begin{cases}
     \approx\omega_{N-2}\ln(1/\epsilon) \text{ if } p=\frac{N+1}{2} \\
      O(1) \text{ if } p>\frac{N+1}{2}
     \end{cases}
$$
and if $p<\frac{N+1}{2}$,
\begin{align*}
& \int_{|y|\le a} \frac{y_i^2y_j^2\, dy} {(\epsilon^2+|y|^2)^{\frac{p(N-1)}{2(p-1)}+1}} \\
= & \epsilon^{1-\frac{N-p}{p-1}} \int_{|y|\le a/\epsilon} \frac{y_i^2y_j^2\, dy}
{(1+|y|^2)^{\frac{p(N-1)}{2(p-1)}+1}} \\
= & 4\omega_{N-4} \int_0^\infty \frac{r^{N-4}\, dr}
{(1+r^2)^{\frac{p(N-1)}{2(p-1)}-2}} \int_0^\infty \frac{y_i^2\, dy_i}
{(1+y_i^2)^{\frac{p(N-1)}{2(p-1)}-\frac{1}{2}}} \int_0^\infty \frac{y_j^2\, dy_j}
{(1+y_j^2)^{\frac{p(N-1)}{2(p-1)}+1}} \\
& + O(1)
\end{align*}
Hence
\begin{align*}
& \int_{|y|\le a} \frac{y_i^2y_j^2\, dy}{(\epsilon^2+|y|^2)^{\frac{p(N-1)}{2(p-1)}+1}} \\
& = \begin{cases} \epsilon^{1-\frac{N-p}{p-1}} \frac{\omega_{N-4}}{2}
\frac{\Gamma\left(\frac{N-3}{2}\right)\Gamma\left(\frac{3}{2}\right)^2
\Gamma\left(\frac{N-2p+1}{2(p-1)}\right)} {\Gamma\left(\frac{p(N-1)}{2(p-1)}+1\right) }
+ O(1) \text{ if } p<\frac{N+1}{2} \\
\approx \omega_{N-2}\ln(1/\epsilon) \text{ if } p=\frac{N+1}{2} \\
O(1) \text{ if } p>\frac{N+1}{2}
\end{cases}
\end{align*}
Once again,
\begin{align*}
\int_{|y|\le a} \frac{|y|^5\, dy}{(\epsilon^2+|y|^2)^{\frac{p(N-1)}{2(p-1)}+1}}
& = \epsilon^{2-\frac{N-p}{p-1}} \omega_{N-2} \int_0^{a/\epsilon} \frac{r^{N+3}\,dr}
{(1+r^2)^{\frac{p(N-1)}{2(p-1)}+1}} \\
& = \begin{cases}
   O(\epsilon^{2-\frac{N-p}{p-1}}) \text{ if } p<\frac{N+2}{3} \\
   O(\ln(1/\epsilon)) \text{ if } p=\frac{N+2}{3} \\
   O(1) \text{ if } p>\frac{N+2}{3}
     \end{cases}
\end{align*}
Using the fact that $\Gamma(\frac32)=\frac{\sqrt{\pi}}{2}$,
$\Gamma(\frac52)=\frac{3\sqrt{\pi}}{4}$, and
$$
\omega_{N-3} = \frac{1}{\sqrt{\pi}}
\frac{\Gamma\left(\frac{N-1}{2}\right)}{\Gamma\left(\frac{N-2}{2}\right)}\omega_{N-2},
\qquad \omega_{N-4} = \frac{1}{\pi}
\frac{\Gamma\left(\frac{N-1}{2}\right)}{\Gamma\left(\frac{N-3}{2}\right)}\omega_{N-2},
$$
we eventually get that
\begin{equation}\label{I4}
I_4 =
 \begin{cases}
  \frac{\omega_{N-2}}{16}\epsilon^{1-\frac{N-p}{p-1}}
 \frac{\Gamma\left(\frac{N-2p+1}{2(p-1)}\right)\Gamma\left(\frac{N-1}{2}\right)}
 {\Gamma\left(\frac{p(N-1)}{2(p-1)}+1\right)}
 \left( \frac{3}{2}\sum \lambda_i^2 + \sum_{i<j} \lambda_i\lambda_j \right) \\
 \hspace{1cm} + \begin{cases}
    O(\epsilon^{2-\frac{N-p}{p-1}}) \text{ if } p<\frac{N+2}{3} \\
    O(\ln(1/\epsilon)) \text{ if } p=\frac{N+2}{3} \\
    O(1) \text{ if } \frac{N+2}{3}<p<\frac{N+1}{2}
     \end{cases}  \\
  \frac{\omega_{N-2}}{2}\ln(1/\epsilon)\left(\frac{1}{2}\sum \lambda_i^2 +
  \sum_{i<j} \lambda_i\lambda_j + o(1)\right) \text{ if } p=\frac{N+1}{2} \\
    O(1) \text{ if } p>\frac{N+1}{2}
     \end{cases}
\end{equation}
We thus obtain
$$
I_2 =
 \begin{cases}
  \epsilon^{1-\frac{N-p}{p-1}} \frac{H(0)\omega_{N-2}}{4}
  \frac{\Gamma\left(\frac{N+1}{2}\right)\Gamma\left(\frac{N-2p+1}{2(p-1)}\right)}
  {\Gamma\left(\frac{p(N-1)}{2(p-1)}\right) } \\
  \hspace{1cm} - \epsilon^{2-\frac{N-p}{p-1}} \frac{\omega_{N-2}}{16}
  \frac{\Gamma\left(\frac{N-1}{2}\right)\Gamma\left(\frac{N-2p+1}{2(p-1)}\right)}
  {\Gamma\left(\frac{p(N-1)}{2(p-1)}\right)}
  \left(\frac{3}{2}\sum \lambda_i^2 + \sum_{i<j} \lambda_i\lambda_j \right) \\
  \hspace{1cm} + \begin{cases}
                  O(\epsilon^{3-\frac{N-p}{p-1}}) \text{ if } p<\frac{N+3}{4} \\
                  O(\ln(1/\epsilon)) \text{ if } p=\frac{N+3}{4} \\
                  O(1) \text{ if } \frac{N+3}{4}<p<\frac{N+1}{2}
                 \end{cases}  \\
  \frac{H(0)\omega_{N-2}}{2}\ln(1/\epsilon)(1 + o(1)) \text{ if } p=\frac{N+1}{2} \\
  O(1) \text{ if } p>\frac{N+1}{2}
 \end{cases}
$$
So the proof of \eqref{NormeGradient} is completed.

To prove \eqref{NormeLp}, we first observe that
\begin{equation}\label{NormeLp_4}
\begin{aligned}
\int_\Omega & h(x) |u_\epsilon|^p\, dx = h(0)\int_\Omega |u_\epsilon|^p\, dx +
O\left(\int_\Omega |x||u_\epsilon|^p\, dx \right) \\
& = h(0)\int_{Q_a} |u_\epsilon|^p\, dx + O\left(\int_{Q_a\setminus\Omega}
|u_\epsilon|^p\, dx + \int_{Q_a} |x||u_\epsilon|^p\, dx \right),
\end{aligned}
\end{equation}
where, as before, $Q_a = \{(y,t)\ |\ |y|\le a \text{ and } 0\le t\le a\}$.

Now,
\begin{align*}
\int_{Q_a} |u_\epsilon|^p dx & = \int_{|y|\le a,0<t\le a} \frac{dydt}
{[(t+\epsilon)^2+|y|^2]^\frac{p(N-p)}{2(p-1)}} + O(1) \\
& = \epsilon^{-\frac{N-p^2}{p-1}} \int_{|y|\le a/\epsilon,0<t\le a/\epsilon} \frac{dydt}
{[(1+t)^2+|y|^2]^\frac{p(N-p)}{2(p-1)}} + O(1) \\
& = \begin{cases}
      O(\ln(1/\epsilon)) \text{ if } p^2=N \\
    O(1) \text{ if } p^2>N
     \end{cases}
\end{align*}
If $p^2<N$, using the change of variable $y=(1+t)z$ and then passing to polar
coordinates, we get
\begin{align*}
\int_{Q_a} |u_\epsilon|^p dx
 & = \epsilon^{-\frac{N-p^2}{p-1}} \omega_{N-2} \int_0^\infty \frac{dt}{(1+t)^{\frac{N-p^2}{p-1}+1}}
 \int_0^\infty \frac{r^{N-2}\, dr}{(1+r^2)^\frac{p(N-p)}{2(p-1)}} + O(1)
\end{align*}
Hence
\begin{equation}\label{NormeLp_1}
\begin{aligned}
\int_{Q_a} |u_\epsilon|^p dx & = \begin{cases} \epsilon^{-\frac{N-p^2}{p-1}}
\frac{p-1}{N-p^2} \omega_{N-2} \frac{\Gamma\left(\frac{N-1}{2}\right)
\Gamma\left(\frac{N-p^2+p-1}{2(p-1)}\right)}
{2\Gamma\left(\frac{p(N-p)}{2(p-1)}\right)} + O(1) \text{ if } p^2<N \\
O(\ln(1/\epsilon)) \text{ if } p^2=N \\
O(1) \text{ if } p^2>N
\end{cases}
\end{aligned}
\end{equation}

On the other hand, using Taylor's formula,
\begin{equation}\label{NormeLp_2}
\begin{aligned}
\int_{Q_a\setminus \Omega} |u_\epsilon|^p dx
& = \int_{|y|\le a} \int_0^{\rho(y)} \frac{dt}{[(t+\epsilon)^2+|y|^2]^\frac{p(N-p)}{2(p-1)}}\, dy + O(1) \\
& = O\left( \int_{|y|\le a} \frac{|y|^2\, dy} {(\epsilon^2+|y|^2)^\frac{p(N-p)}{2(p-1)}}\, dy \right) + O(1) \\
& = \epsilon^{1-\frac{N-p^2}{p-1}}  O\left(\int_0^{a/\epsilon} \frac{r^N\, dr}{(1+r^2)^\frac{p(N-p)}
{2(p-1)}}\right) + O(1) \\
& = \begin{cases}
      O(\epsilon^{1-\frac{N-p^2}{p-1}}) \text{ if } p<\frac{-1+\sqrt{4N+5}}{2} \\
      O(\ln(1/\epsilon)) \text{ if } p=\frac{-1+\sqrt{4N+5}}{2} \\
      O(1) \text{ if } p>\frac{-1+\sqrt{4N+5}}{2}
     \end{cases}
\end{aligned}
\end{equation}
Similarly,
\begin{equation}\label{NormeLp_3}
\begin{aligned}
\int_{Q_a} |x| |u_\epsilon|^p dx
& = \int_{Q_a} \frac{|(y,t)|}{[(t+\epsilon)^2+|y|^2]^\frac{p(N-p)}{2(p-1)}}\, dydt + O(1) \\
& = \epsilon^{1-\frac{N-p^2}{p-1}} \int_{Q_{a/\epsilon}}
\frac{|(y,t)|}{[(1+t)^2+|y|^2]^\frac{p(N-p)}{2(p-1)}}\, dydt + O(1) \\
& = \begin{cases}
       O(\epsilon^{1-\frac{N-p^2}{p-1}}) \text{ if } p<\frac{-1+\sqrt{4N+5}}{2} \\
       O(\ln(1/\epsilon)) \text{ if } p=\frac{-1+\sqrt{4N+5}}{2} \\
       O(1) \text{ if } p>\frac{-1+\sqrt{4N+5}}{2}
     \end{cases}
\end{aligned}
\end{equation}

Combining \eqref{NormeLp_4}, \eqref{NormeLp_1}, \eqref{NormeLp_2} and
\eqref{NormeLp_3}, gives \eqref{NormeLp}.

Finally, to prove \eqref{NormeCritique}, we first observe that
$$
\int_{\p\Omega} |u_\epsilon|^{p_*}\, dS = \int_{Q_a} |u_\epsilon|^{p_*}\, dS
$$
for small $\epsilon$ and so
\begin{equation*}
\begin{aligned}
\int_{\p\Omega} |u_\epsilon|^{p_*}\, dS = & \int_{|y|\le a}
\frac{\sqrt{1+|\nabla\rho|^2}}{[(\epsilon+\rho(y))^2+|y|^2]^\frac{p(N-1)}{2(p-1)}}\, dy \\
= & \int_{|y|\le a} \frac{1+\frac12 |\nabla\rho|^2 + O(|y|^4)}
{(\epsilon^2+|y|^2)^\frac{p(N-1)}{2(p-1)}} \Big[1-\frac{p(N-1)}{2(p-1)}
\frac{\rho(2\epsilon+\rho)}{\epsilon^2+|y|^2} \\
& - c_{N,p}\frac{\rho^2(2\epsilon+\rho)^2}{(\epsilon^2+|y|^2)^2}  +
O\left(\frac{\rho^3(2\epsilon+\rho)^3}{(\epsilon^2+|y|^2)^3}\right)\Big]\, dy,
\end{aligned}
\end{equation*}
where
$$
c_{N,p} = -\frac{p(N-1)}{4(p-1)}\left[ \frac{p(N-1)}{2(p-1)}+1\right].
$$

Hence
\begin{align*}
& \int_{\p\Omega} |u_\epsilon|^{p_*}\, dS = \\
= & \int_{|y|\le a} \frac{dy}{(\epsilon^2+|y|^2)^\frac{p(N-1)}{2(p-1)}}\, dy
- \epsilon^\frac{p(N-1)}{p-1} \int_{|y|\le a} \frac{\rho(y)\, dy}
{(\epsilon^2+|y|^2)^{1+\frac{p(N-1)}{2(p-1)}}} \\
& + \frac12 \int_{|y|\le a} \frac{|\nabla\rho|^2\, dy}
{(\epsilon^2+|y|^2)^\frac{p(N-1)}{2(p-1)}} - \frac{p(N-1)}{2(p-1)}
\int_{|y|\le a} \frac{\rho^2(y)\, dy}{(\epsilon^2+|y|^2)^{1+\frac{p(N-1)}{2(p-1)}}} \\
& - 4\epsilon^2 c_{N,p} \int_{|y|\le a} \frac{\rho^2(y)\, dy}
{(\epsilon^2+|y|^2)^{2+\frac{p(N-1)}{2(p-1)}}}  \\
& + O\left(\int_{|y|\le a} \frac{|y|^4\, dy}
{(\epsilon^2+|y|^2)^\frac{p(N-1)}{2(p-1)}}\, dy + \epsilon \int_{|y|\le a}
\frac{|y|^4\, dy}{(\epsilon^2+|y|^2)^{1+\frac{p(N-1)}{2(p-1)}}}\, dy \right)  \\
& = I_5 - \epsilon^\frac{p(N-1)}{p-1} I_7 + \frac12 I_6 - \frac{p(N-1)}{2(p-1)}
I_8 - 4\epsilon^2 c_{N,p} I_9 + O(I_{10}).
\end{align*}

We first compute $I_5$ as follows:
\begin{equation}\label{I5}
\begin{aligned}
I_5 & = \int_{|y|\le a} \frac{dy}{(\epsilon^2+|y|^2)^\frac{p(N-1)}{2(p-1)}} =
\omega_{N-2} \epsilon^{-1-\frac{N-p}{p-1}} \int_0^{a/\epsilon} \frac{r^{N-2}\, dr}
{(1+r^2)^\frac{p(N-1)}{2(p-1)}} \\
& = \omega_{N-2} \epsilon^{-1-\frac{N-p}{p-1}} \int_0^\infty \frac{r^{N-2}\, dr}
{(1+r^2)^\frac{p(N-1)}{2(p-1)}} + O(1) \\
& = \omega_{N-2} \epsilon^{-1-\frac{N-p}{p-1}}
\frac{\Gamma\left(\frac{N-1}{2}\right) \Gamma\left(\frac{N-1}{2(p-1)}\right)}
{2\Gamma\left(\frac{p(N-1)}{2(p-1)}\right)} + O(1).
\end{aligned}
\end{equation}
According to \eqref{E1} and \eqref{E2}, using the relation
$\Gamma\left(\frac{N+1}{2}\right)=\frac{N-1}{2}\Gamma\left(\frac{N-1}{2}\right)$,
we have
\begin{equation}\label{I6}
\begin{aligned}
I_6 & = \int_{|y|\le a} \frac{|\nabla \rho|^2\, dy}{(\epsilon^2+|y|^2)^\frac{p(N-1)}{2(p-1)}} \\
& = \sum \lambda_i^2 \int_{|y|\le a} \frac{|y_i|^2\, dy}
{(\epsilon^2+|y|^2)^\frac{p(N-1)}{2(p-1)}} + O\left(\int_{|y|\le a} \frac{|y|^4\, dx}
{(\epsilon^2+|y|^2)^\frac{p(N-1)}{2(p-1)}} \right) \\
& = \frac{\sum \lambda_i^2}{N-1} \int_{|y|\le a} \frac{|y|^2\, dy}
{(\epsilon^2+|y|^2)^\frac{p(N-1)}{2(p-1)}} + O\left(\int_{|y|\le a} \frac{|y|^4\, dx}
{(\epsilon^2+|y|^2)^\frac{p(N-1)}{2(p-1)}} \right) \\
& = \begin{cases} \frac14 \sum \lambda_i^2 \omega_{N-2}
\epsilon^{1-\frac{N-p}{p-1}}
\frac{\Gamma\left(\frac{N-1}{2}\right)\Gamma\left(\frac{N-2p+1}{2(p-1)}\right)}
{\Gamma\left(\frac{p(N-1)}{2(p-1)}\right)} + \begin{cases}
O(\epsilon^{3-\frac{N-p}{p-1}}) \text{ if } p<\frac{N+3}{4} \\
O(\ln(1/\epsilon)) \text{ if } p=\frac{N+3}{4} \\
O(1) \text{ if } \frac{N+1}{2}>p>\frac{N+3}{4}
\end{cases} \\
\frac{\omega_{N-2}\sum \lambda_i^2}{N-1} \ln(1/\epsilon)) \text{ if } p=\frac{N+1}{2} \\
O(1) \text{ if } p>\frac{N+1}{2}
\end{cases}
\end{aligned}
\end{equation}

By radial symmetry, we have
\begin{align*}
I_7 = & \int_{|y|\le a} \frac{\rho(y)\, dy}{(\epsilon^2+|y|^2)^{1+\frac{p(N-1)}{2(p-1)}}} \\
= & \frac{\sum \lambda_i}{2(N-1)} \int_{|y|\le a} \frac{|y|^2\, dy}
{(\epsilon^2+|y|^2)^{1+\frac{p(N-1)}{2(p-1)}}} + O\left(\int_{|y|\le a}
\frac{|y|^4\, dy} {(\epsilon^2+|y|^2)^{1+\frac{p(N-1)}{2(p-1)}}} \right) \\
= & \frac{\omega_{N-2} \sum \lambda_i}{2(N-1)} \epsilon^{-1-\frac{N-p}{p-1}}
\int_0^{a/\epsilon} \frac{r^N\, dr}{(1+r^2)^{1+\frac{p(N-1)}{2(p-1)}}} \\
& + \epsilon^{-\frac{N-p}{p-1}} O\left(\int_0^{a/\epsilon} \frac{r^{N+2}\, dr}
{(1+r^2)^{1+\frac{p(N-1)}{2(p-1)}}} \right)\\
= & \frac{\omega_{N-2} \sum \lambda_i}{2(N-1)} \epsilon^{-1-\frac{N-p}{p-1}}
\int_0^\infty \frac{r^N\, dr} {(1+r^2)^{1+\frac{p(N-1)}{2(p-1)}}} +
\begin{cases}
O(\epsilon^{1-\frac{N-p}{p-1}} ) \text{ if } p<\frac{N+1}{2} \\
O(\ln(1/\epsilon)) \text{ if } p=\frac{N+1}{2} \\
O(\epsilon^{-\frac{N-p}{p-1}}) \text{ if } p>\frac{N+1}{2}
\end{cases}
\end{align*}
and so
\begin{equation}\label{I7}
\begin{aligned}
I_7 = & \frac{\omega_{N-2} \sum \lambda_i}{8} \epsilon^{-1-\frac{N-p}{p-1}}
\frac{\Gamma\left(\frac{N-1}{2}\right)\Gamma\left(\frac{N-1}{2(p-1)}\right)}
{\Gamma\left(1+\frac{p(N-1)}{2(p-1)}\right)} \\
& + \begin{cases}
O(\epsilon^{1-\frac{N-p}{p-1}}) \text{ if } p<\frac{N+1}{2} \\
O(\ln(1/\epsilon)) \text{ if } p=\frac{N+1}{2} \\
O(\epsilon^{-\frac{N-p}{p-1}}) \text{ if } p>\frac{N+1}{2}
\end{cases}
\end{aligned}
\end{equation}
To compute $I_9$ we proceed as in the computations of $I_4$, i.e.
\begin{align*}
I_9 = & \int_{|y|\le a} \frac{\rho^2(y)\, dy} {(\epsilon^2+|y|^2)^{2+\frac{p(N-1)}{2(p-1)}} } \\
= & \frac14 \sum \lambda_i^2 \int_{|y|\le a} \frac{y_1^4\, dy}
{(\epsilon^2+|y|^2)^{2+\frac{p(N-1)}{2(p-1)}}}\\
& + \frac12 \sum_{i<j} \lambda_i\lambda_j \int_{|y|\le a} \frac{y_i^2 y_j^2\,
dy} {(\epsilon^2+|y|^2)^{2+\frac{p(N-1)}{2(p-1)}}} + O\left( \int_{|y|\le a}
\frac{|y|^5\, dy} {(\epsilon^2+|y|^2)^{2+\frac{p(N-1)}{2(p-1)}} } \right).
\end{align*}
Now
\begin{align*}
\int_{|y|\le a} & \frac{y_1^4\, dy}
{(\epsilon^2+|y|^2)^{2+\frac{p(N-1)}{2(p-1)}}} = \epsilon^{-\frac{N-1}{p-1}}
\int_{\R^{N-1}} \frac{y_1^4\, dy}
{(1+|y|^2)^{2+\frac{p(N-1)}{2(p-1)}}} + O(1) \\
& = 2 \epsilon^{-\frac{N-1}{p-1}} \omega_{N-3} \int_0^\infty \frac{r^{N-3}\,
dr} {(1+r^2)^{\frac{p(N-1)}{2(p-1)}-\frac12}} \int_0^\infty \frac{s^4\, ds}
{(1+s^2)^{2+\frac{p(N-1)}{2(p-1)}}} + O(1) \\
& = \frac{3\omega_{N-2}}{8} \epsilon^{-\frac{N-1}{p-1}}
\frac{\Gamma\left(\frac{N-1}{2}\right) \Gamma\left(\frac{N-1}{2(p-1)}\right)}
{\Gamma\left(2+\frac{p(N-1)}{2(p-1)}\right)} + O(1),
\end{align*}

\begin{align*}
\int_{|y|\le a} & \frac{y_i^2 y_j^2\, dy}
{(\epsilon^2+|y|^2)^{2+\frac{p(N-1)}{2(p-1)}} } = \epsilon^{-\frac{N-1}{p-1}}
\int_{\R^{N-1}} \frac{y_i^2 y_j^2\, dy}{(1+|y|^2)^{2+\frac{p(N-1)}{2(p-1)}}} + O(1) \\
= & 4\epsilon^{-\frac{N-1}{p-1}} \omega_{N-4} \int_0^\infty \frac{r^{N-4}\, dr}
{(1+r^2)^{\frac{p(N-1)}{2(p-1)}-1}} \int_0^\infty \frac{y_i^2\, dy_i}
{(1+y_i^2)^{\frac12+\frac{p(N-1)}{2(p-1)}} } \\
& \times \int_0^\infty \frac{y_j^2 dy_j}{(1+y_j^2)^{2+\frac{p(N-1)}{2(p-1)}} } + O(1) \\
= & \frac{\omega_{N-2}}{8} \epsilon^{-\frac{N-1}{p-1}}
\frac{\Gamma\left(\frac{N-1}{2}\right)\Gamma\left(\frac{N-1}{2(p-1)}\right)}
{\Gamma\left(2+\frac{p(N-1)}{2(p-1)}\right)} + O(1),
\end{align*}
and
\begin{align*}
\int_{|y|\le a} \frac{|y|^5\, dy}
{(\epsilon^2+|y|^2)^{2+\frac{p(N-1)}{2(p-1)}}}
& = \epsilon^{-\frac{N-p}{p-1}} \omega_{N-2}  \int_0^{a/\epsilon} \frac{r^{N+3}\, dr}
{(1+r^2)^{2+\frac{p(N-1)}{2(p-1)}}} \\
& = O(\epsilon^{-\frac{N-p}{p-1}})
\end{align*}
Hence
\begin{equation}\label{I9}
\begin{aligned}
I_9 = & \frac{\omega_{N-2}}{16} \epsilon^{-\frac{N-1}{p-1}}
\frac{\Gamma\left(\frac{N-1}{2}\right) \Gamma\left(\frac{N-1}{2(p-1)}\right)}
{\Gamma\left(2+\frac{p(N-1)}{2(p-1)}\right)} \left(\frac32 \sum \lambda_i^2 +
\sum_{i<j} \lambda_i\lambda_j \right)\\
& + O(\epsilon^{-\frac{N-p}{p-1}}).
\end{aligned}
\end{equation}

Finally, for $I_{10}$ we have,
\begin{align*}
I_{10} = & \epsilon^{3-\frac{N-p}{p-1}} \omega_{N-2} \int_0^{a/\epsilon}
\frac{r^{N+2}\, dr}{(1+r^2)^\frac{p(N-1)}{2(p-1)}} + \epsilon^{2-\frac{N-p}{p-1}}
\omega_{N-2} \int_0^{a/\epsilon} \frac{r^{N+2}\, dr}{(1+r^2)^{1+\frac{p(N-1)}{2(p-1)}}} \\
= & \begin{cases}
       O(\epsilon^{3-\frac{N-p}{p-1}}) \text{ if } p<\frac{N+3}{4} \\
       O(\ln(1/\epsilon)) \text{ if } p=\frac{N+3}{4} \\
       O(1) \text{ if } p>\frac{N+3}{4}
     \end{cases}
+ \begin{cases}
       O(\epsilon^{2-\frac{N-p}{p-1}}) \text{ if } p<\frac{N+1}{2} \\
       O(\epsilon\ln(1/\epsilon)) \text{ if } p=\frac{N+1}{2} \\
       O(\epsilon) \text{ if } p>\frac{N+1}{2}
     \end{cases}
\end{align*}
and so
\begin{equation}\label{I10}
I_{10} =  \begin{cases}
          O(\epsilon^{2-\frac{N-p}{p-1}}) \text{ if } p\le\frac{N+2}{3} \\
          O(1) \text{ if } p>\frac{N+2}{3}
     \end{cases}
\end{equation}

Putting these estimates together, we arrive at \eqref{NormeCritique}. This
completes the proof of Step \ref{step_calcul1}.
\end{proof}

\begin{Step}\label{Step3}
We have, for any dimension $N\ge 2$,
$$
K_p^{-1}\frac{\displaystyle \int_{\Omega} |\nabla u_\epsilon|^p +
|u_\epsilon|^p\, dx}{\displaystyle \Big(\int_{\p\Omega} |u_\epsilon|^{p_*}\,
dS\Big)^{p/p_*}} =
\begin{cases}
        1+O(\epsilon^\frac{N-p}{p-1}) \text{ if } p>\frac{N+1}{2} \\
        1 - \frac{N-1}{2}H(0)\epsilon \ln(1/\epsilon) + o(\epsilon \ln(1/\epsilon))
        \text{ if } p=\frac{N+1}{2}
\end{cases}
$$
and, if $p<\frac{N+1}{2}$, for dimension $N=2,3,4$
\begin{align*}
K_p^{-1}\frac{\displaystyle \int_{\Omega} |\nabla u_\epsilon|^p +
|u_\epsilon|^p\, dx}{\displaystyle \Big(\int_{\p\Omega} |u_\epsilon|^{p_*}\,
dS\Big)^{p/p_*}} =
& 1 -\frac{(N-p)(p-1)}{N-2p+1} H(0) \epsilon  \\
& + \begin{cases}
  \frac{D}{A_1}\epsilon^p + \begin{cases}
                                E \epsilon^2 + O(\epsilon^{1+p}) \text{ if } p<\frac{N+2}{3} \\
                                O(\epsilon^\frac{N-p}{p-1}) \text{ if } \frac{N+2}{3}\le p<\sqrt{N}
                              \end{cases} \\
     O(\epsilon^\frac{N-p}{p-1}\ln(1/\epsilon))  \text{ if } p=\sqrt{N} \\
     O(\epsilon^\frac{N-p}{p-1}) \text{ if } \sqrt{N}<p<\frac{N+1}{2}
   \end{cases}
\end{align*}
where
$$
E = \frac{(N-p)(p-1)}{4(N-1)(N-2p+1)} \left\{\frac{p+N-2}{N-1}\sum \lambda_i^2
- 2\sum_{i<j}\lambda_i\lambda_j \right\}.
$$
Also, for dimensions $N\ge 5$,
\begin{align*}
K_p^{-1}\frac{\displaystyle \int_{\Omega} |\nabla u_\epsilon|^p +
|u_\epsilon|^p\, dx}{\displaystyle \Big(\int_{\p\Omega} |u_\epsilon|^{p_*}\,
dS\Big)^{p/p_*}} = & 1 -\frac{(N-p)(p-1)}{N-2p+1} H(0) \epsilon  \\
& + \begin{cases}
       E\epsilon^2 + \begin{cases}
                      \frac{D}{A_1}\epsilon^p + \begin{cases}
                                                  o(\epsilon^2) \text{ if } p\le 2 \\
                                                  o(\epsilon^p) \text{ if } 2\le p <\sqrt{N}
                                                 \end{cases} \\
                       o(\epsilon^2) \text{ if } \sqrt{N}\le p<\frac{N+2}{3}
                      \end{cases} \\
     O(\epsilon^2) \text{ if } \frac{N+2}{3}\le p<\frac{N+1}{2}
   \end{cases}
\end{align*}
\end{Step}

\begin{proof}[Proof of Step \ref{Step3}]
Noting that
$$
\frac{A_1}{B_1^\frac{N-p}{N-1}} = K_p^{-1},
$$
we have, when e.g. $n\ge 6$ and $p\le 2$, that
\begin{align*}
K_p^{-1}&\frac{\displaystyle \int_{\Omega} |\nabla u_\epsilon|^p +
|u_\epsilon|^p\, dx}{\displaystyle \Big(\int_{\p\Omega} |u_\epsilon|^{p_*}\,
dS\Big)^{p/p_*}} = 1+\left(\frac{A_2}{A_1} - \frac{N-p}{N-1}
\frac{B_2}{B_1}\right)\epsilon
+ \frac{D}{A_1}\epsilon^p \\
& + \left\{ \frac{N-p}{N-1} \left[\frac12 \left(\frac{N-p}{N-1}+1\right)
\left(\frac{B_2}{B_1}\right)^2 - \frac{B_3}{B_1} -
\frac{B_2}{B_1}\frac{A_2}{A_1}\right] + \frac{A_3}{A_1} \right\} \epsilon^2 +
o(\epsilon^2).
\end{align*}
Using the fact that
\begin{align*}
& \Gamma\left(\frac{N+1}{2}\right) = \Gamma\left(\frac{N-1}{2}+1\right) =
\frac{N-1}{2}\Gamma\left(\frac{N-1}{2}\right)\\
& \Gamma\left(\frac{N-1}{2(p-1)}\right) =
\Gamma\left(\frac{N-2p+1}{2(p-1)}+1\right) =
\frac{N-2p+1}{2(p-1)}\Gamma\left(\frac{N-2p+1}{2(p-1)} \right),
\end{align*}
we get
\begin{align*}
& \frac{A_2}{A_1} = -\frac12 \frac{N-p}{N-2p+1} \sum \lambda_i, \\
& \frac{A_3}{A_1} = \frac14 \frac{N-p}{N-2p+1} \left\{\frac32 \sum \lambda_i^2 - 2
\sum_{i<j} \lambda_i\lambda_j \right\}, \\
& \frac{B_2}{B_1} = -\frac12 \sum \lambda_i, \\
& \frac{B_3}{B_1} = \frac{1}{8(N-2p+1)} \left\{(3N-5p+2)\sum \lambda_i^2 - 4(N-p)
\sum_{i<j} \lambda_i\lambda_j \right\}, \\
& \frac{D}{A_1} = \begin{cases}
                    \frac{2h(0)}{(N-3)(N-4)} \text{ if } p=2 \\
                    \text{has same sign as } \frac{h(0)}{N-p^2} \text{ otherwise}.
                   \end{cases}
\end{align*}
Hence
$$
\frac{A_2}{A_1}-\frac{N-p}{N-1}\frac{B_2}{B_1} = -\frac{(N-p)(p-1)}{N-2p+1}
H(0)
$$
and
\begin{align*}
\frac{N-p}{N-1} &\left[\frac12 \left(\frac{N-p}{N-1} + 1\right)
\left(\frac{B_2}{B_1}\right)^2 - \frac{B_3}{B_1} - \frac{B_2}{B_1}
\frac{A_2}{A_1}\right] + \frac{A_3}{A_1} \\
= & \frac{(N-p)(p-1)}{4(N-1)(N-2p+1)} \left\{\frac{p+N-2}{N-1} \sum \lambda_i^2
- 2\sum_{i<j}\lambda_i\lambda_j \right\},
\end{align*}
which gives the result. We get the others equalities in much the same way.
\end{proof}

\begin{proof}[\bf Proof of Theorem \ref{thm1}]
At this point is just a combination of Steps 1 and 2.
\end{proof}

\end{document}